\documentclass[11pt, leqno]{article}
\usepackage{amsmath, amsfonts}
\begin{document}

\author{Ajai Choudhry and Arman Shamsi Zargar}
\title{A parametrised family of Mordell curves}
\date{}
\maketitle

\begin{abstract}
An elliptic curve defined by an equation of the type $y^2=x^3+d$ is called a Mordell curve. We obtain a parametrised family of Mordell curves whose rank, in general, is at least three, and whose torsion group is $\mathbb{Z}/3\mathbb{Z}$.
\end{abstract}



Ever since Fermat's assertion that the only solution in positive integers of the equation $y^2=x^3-2$ is $(x,\,y)=(3,\,5)$ \cite{Mo2}, the diophantine equation,
\begin{equation}
y^2=x^3+d, \label{Moeq}
\end{equation}
has been subjected to extensive investigations (\cite{Ben},  \cite{Ell},  \cite{Geb}, \cite{Ha}, \cite{LJB}, \cite{Lj}, \cite{Mo1},  \cite{Mo2a},  \cite[Chapter 26, pp. 238--254]{Mo3}, \cite{Po}, \cite{Yo}). Eq.~\eqref{Moeq} has now been solved for all integer values of $d$ with  $|d| \leq 10^7$ \cite{Ben}.

 The elliptic curve represented by Eq.~\eqref{Moeq} is known as a Mordell curve. Further, for various integer values of $d$, we know the structure of the torsion subgroup of the group of rational points  on the  Mordell curve \cite[Theorem 5.3, p. 134]{Kna}. It is noteworthy that whenever the integer $d$ is a nonzero perfect square different from 1, the torsion subgroup is necessarily  $\mathbb{Z}/3\mathbb{Z}$.

In this paper, we construct  a parametrised family  of  Mordell curves defined by an equation of the type,
\begin{equation}
y^2=x^3+k^2. \label{Moeqksq}
\end{equation}
We shall show that  the  rank of the elliptic curves belonging to this family is, in general,  at least three. The torsion subgroup of all curves defined by the equation \eqref{Moeqksq} is $\mathbb{Z}/3\mathbb{Z}$ with $(0,\,k)$ being a torsion  point of order 3.  Despite the considerable literature related to Eq.~\eqref{Moeq}, it appears that such a  parametrised family of curves has not been obtained earlier.

We will first solve the system of diophantine equations,
\begin{align}
v_1^2&=u_1^3+k^2, \label{eq1}\vspace{.2cm}\\
v_2^2&=u_2^3+k^2, \label{eq2}\vspace{.2cm}\\
v_3^2&=u_3^3+k^2. \label{eq3}
\end{align}
On writing,
\begin{equation}
u_1=am,\quad u_2=bm, \quad u_3=cm, \label{valu}
\end{equation}
and eliminating $m$ first between Eq.~\eqref{eq1} and Eq.~\eqref{eq2}, and then  between Eq.~\eqref{eq1} and Eq.~\eqref{eq3}, we get the following two equations:
\begin{align}
b^3(v_1^2-k^2)&=a^3(v_2^2-k^2), \label{eq4} \\
c^3(v_1^2-k^2)&=a^3(v_3^2-k^2). \label{eq5}
\end{align}
To solve equations \eqref{eq4} and \eqref{eq5}, we write,
\begin{equation*}
v_1=w_1t+k,\quad v_2=w_2t+k,\quad v_3=w_3t+k,
\end{equation*}
when each of the two equations \eqref{eq4} and \eqref{eq5} can be readily solved to get a nonzero solution for $t$. Equating these two values of $t$, we get the condition,
\begin{equation}
(b^3w_1-a^3w_2)(c^3w_1^2-a^3w_3^2) = (c^3w_1-a^3w_3)(b^3w_1^2-a^3w_2^2). \label{condw}
\end{equation}

Now Eq.~\eqref{condw} is a homogeneous cubic equation in the variables $w_1,\,w_2$ and $w_3$ and it represents a cubic curve in the projective plane. Further,   a rational point on this curve is easily seen to be $(w_1,\,w_2,\,w_3)=(a^3,\,b^3,\,c^3).$ The tangent to the cubic curve \eqref{condw} at the aforementioned rational point necessarily intersects the curve \eqref{condw} at another rational point which is thus easily found, and is given by,
\begin{equation*}
w_1 = a^3(b^3+c^3-a^3),\;\; w_2 = b^3(c^3+a^3-b^3),\;\; w_3 = c^3(a^3+b^3-c^3). 
\end{equation*}

With these values of $w_1,\,w_2,\,w_3$, we obtain the following solution of the simultaneous  equations \eqref{eq4} and  \eqref{eq5}:
\begin{equation}
\begin{aligned}
v_1&=-(3a^6-2a^3b^3-2a^3c^3-b^6+2b^3c^3-c^6)r,\\
v_2&=(a^6+2a^3b^3-2a^3c^3-3b^6+2b^3c^3+c^6)r,\\
v_3&=(a^6-2a^3b^3+2a^3c^3+b^6+2b^3c^3-3c^6)r,\\
k&=(a^6-2a^3b^3-2a^3c^3+b^6-2b^3c^3+c^6)r,
\end{aligned}
\label{valvk}
\end{equation}
where $a,\,b,\,c$ and $r$ are arbitrary parameters. Substituting the values of $v_1,\,v_2,\,v_3$ and $k$ given by \eqref{valvk} and the value of $u_1$ given by \eqref{valu}  in \eqref{eq1}, we get,
\begin{equation}
a^3m^3=-8a^3r^2(a^3+b^3-c^3)(b^3+c^3-a^3)(c^3+a^3-b^3). \label{condmr}
\end{equation}
Now Eq.~\eqref{condmr} is readily solved by taking $m=r$, when we get,
\begin{equation*}
r=-8(a^3+b^3-c^3)(b^3+c^3-a^3)(c^3+a^3-b^3).
\end{equation*}

We thus obtain a  solution of the system of equations \eqref{eq1}, \eqref{eq2} and \eqref{eq3} which is given by,
\begin{multline}
k = -8(a^3+b^3-c^3)(b^3+c^3-a^3)(c^3+a^3-b^3)\\
\times (a^6-2a^3b^3-2a^3c^3+b^6-2b^3c^3+c^6), \quad \quad \quad \quad \quad \quad \quad \quad \label{valk}
\end{multline}
and
\begin{equation}
\begin{aligned}
u_1&=-8a(a^3+b^3-c^3)(b^3+c^3-a^3)(c^3+a^3-b^3),\\
u_2&=-8b(a^3+b^3-c^3)(b^3+c^3-a^3)(c^3+a^3-b^3),\\
u_3&=-8c(a^3+b^3-c^3)(b^3+c^3-a^3)(c^3+a^3-b^3),\\
v_1&=8(a^3+b^3-c^3)(b^3+c^3-a^3)(c^3+a^3-b^3)\\
& \quad \quad \times (3a^6-2a^3b^3-2a^3c^3-b^6+2b^3c^3-c^6),\\
v_2&=-8(a^3+b^3-c^3)(b^3+c^3-a^3)(c^3+a^3-b^3)\\
& \quad \quad \times (a^6+2a^3b^3-2a^3c^3-3b^6+2b^3c^3+c^6),\\
v_3&=-8(a^3+b^3-c^3)(b^3+c^3-a^3)(c^3+a^3-b^3)\\
& \quad \quad \times (a^6-2a^3b^3+2a^3c^3+b^6+2b^3c^3-3c^6),
\end{aligned}
\label{soleq123}
\end{equation}
where $a,\,b,\,c$ are arbitrary parameters.

It follows that when $k$ is given by \eqref{valk},  there are three rational points $P_1(a,\,b,\,c),$ $P_2(a,\,b,\,c) $ and $P_3(a,\,b,\,c)$ on the elliptic curve \eqref{Moeqksq}  with co-ordinates $(u_i,\,v_i),\;i=1,\,2,\,3,$ where the values of $u_i,\,v_i,i=1,\,2,\,3,$ are given by \eqref{soleq123}.

We will now apply a theorem of Silverman \cite[Theorem 11.4, p. 271]{Sil} to show that these points are linearly independent. For this,  we must find a specialisation $(a,b,c)=(a_0,b_0,c_0)$ such that the points $P_1(a_0,b_0,c_0)$,  $P_2(a_0,b_0,c_0)$, and $P_3(a_0,b_0,c_0)$ are linearly independent on the specialised curve over $\mathbb{Q}$.

We take $(a,b,c)=(1,2,3)$, when we get the elliptic curve,
\begin{equation}
y^2=x^3+28592640^2, \label{Moeqspl}
\end{equation}
on which we get the three points,
\begin{align*}
P_1(1,2,3)&=(97920,41909760), \\
P_2(1,2,3)&=(195840,91261440), \\
P_3(1,2,3)&=(293760,161763840),
\end{align*}
each of which  is of infinite order. The regulator of these three points, as determined by the software SAGE \cite{Sag} is  $33.9574760167017$. As this is nonzero,  it follows from a well-known theorem \cite[Theorem 8.1, p. 242]{SZ} that these three points are linearly independent.  Hence the rank of the  Mordell curve \eqref{Moeqspl} is at least three.

It follows that, in general, the rank of the elliptic curves belonging to  the parametrised family of Mordell curves \eqref{Moeqksq}, where $k$ is given by \eqref{valk}, is at least three.

\noindent Postal Address 1: Ajai Choudhry, 
\newline \hspace{1.05 in} 13/4 A Clay Square,
\newline \hspace{1.05 in} Lucknow - 226001, INDIA.
\newline \noindent  E-mail: ajaic203@yahoo.com

\noindent Postal Address 2: Arman Shamsi Zargar, 
\newline \hspace{1.05 in} Young Researchers and Elite Club, Ardabil Branch,
\newline \hspace{1.05 in} Islamic Azad University, Ardabil, Iran.
\newline \noindent  E-mail: shzargar.arman@gmail.com

\end{document}